\documentclass[12pt]{article} 
\def\R{\hbox{{\rm I}\kern-0.2em{\rm R}\kern0.2em}}%mathematical R for reals
\def\D{\hbox{{\rm I}\kern-0.2em{\rm D}\kern0.2em}}

\def\be{\begin{equation}}
\def\ee{\end{equation}}

\def\({\left(}
\def\){\right)}
\def\[{\left[}
\def\]{\right]}
\def\bc{\begin{center}}
\def\ec{\end{center}}

\begin{document}

{\large \bf Invariants for systems of two linear hyperbolic-type
equations by complex methods}

\textbf{A. Aslam$^{a}$, M. Safdar$^{b}$ and F. M. Mahomed$^{c}$}

$^{a}$School of Natural Sciences\\
National University of Sciences and Technology\\
Campus H-12, 44000, Islamabad, Pakistan

$^{b}$School of Mechanical and Manufacturing Engineering\\
National University of Sciences and Technology\\
Campus H-12, 44000, Islamabad, Pakistan

$^{c}$Differential Equations, Continuum Mechanics and
Applications\\
School of Computational and Applied Mathematics\\
University of the Witwatersrand, Wits 2050, South Africa

adnangrewal@yahoo.com, safdar.camp@gmail.com,
Fazal.Mahomed@wits.ac.za

\textbf{Abstract.} Invariants of general linear system of two
hyperbolic partial differential equations (PDEs) are derived under
transformations of the dependent and independent variables by real
infinitesimal method earlier. Here a subclass of the general system
of linear hyperbolic PDEs is investigated for the associated
invariants, by complex as well as real methods. The complex
procedure relies on the correspondence of systems of PDEs with the
base complex equation. Complex invariants of the base complex PDEs
are shown to reveal invariants of the corresponding systems. A
comparison of all the invariant quantities obtained by complex and
real methods for this class, is presented which shows that the
complex procedure provides a few invariants different from those
extracted by real symmetry analysis.

\section{Introduction}
Invariants of a system of two linear hyperbolic equations
\begin{eqnarray}
u_{tx}+a_{1}(t,x)u_{t}+a_{2}(t,x)v_{t}+b_{1}(t,x)u_{x}+b_{2}(t,x)v_{x}+c_{1}(t,x)u+c_{2}(t,x)v=0,\nonumber\\
v_{tx}+a_{3}(t,x)u_{t}+a_{4}(t,x)v_{t}+b_{3}(t,x)u_{x}+b_{4}(t,x)v_{x}+c_{3}(t,x)u+c_{4}(t,x)v=0,\label{rehypsysgen}
\end{eqnarray}
where $t$ and $x$ in the subscripts denote partial derivatives with
respect to these independent variables, under an invertible change
of the dependent variables have been derived in \cite{hsysp} by real
infinitesimal method. Moreover, joint invariants were also worked
out there, under transformations of both the dependent and
independent variables. The most general group of equivalence
transformations \cite{book1,book2}, i.e. invertible change of the
dependent and independent variables that maps a system of hyperbolic
equations to itself with, in general different coefficients, is
obtained first and then employed to find associated invariants
\cite{hsysp}. Complex symmetry analysis (CSA) is invoked to
investigate the semi-invariants associated with a subclass of the
system of hyperbolic equations (\ref{rehypsysgen}) under a change of
only the dependent variables \cite{comh}. This subclass of systems
is represented by the following two hyperbolic-type equations
\begin{eqnarray}
u_{tx}+\alpha_{1}(t,x)u_{t}-\alpha_{2}(t,x)v_{t}+\beta_{1}(t,x)u_{x}-\beta_{2}(t,x)v_{x}+\gamma_{1}(t,x)u-\gamma_{2}(t,x)v=0,\nonumber\\
v_{tx}+\alpha_{2}(t,x)u_{t}+\alpha_{1}(t,x)v_{t}+\beta_{2}(t,x)u_{x}+\beta_{1}(t,x)v_{x}+\gamma_{2}(t,x)u+\gamma_{1}(t,x)v=0.\label{chypsysgen}
\end{eqnarray}
This system of hyperbolic-type equations corresponds to a scalar
complex hyperbolic equation
\begin{eqnarray}
w_{tx}+\alpha(t,x)w_{t}+\beta(t,x)w_{x}+\gamma(t,x)w=0,\label{scacomhyp}
\end{eqnarray}
if $\alpha_{1}+i\alpha_{2}=\alpha$, $\beta_{1}+i\beta_{2}=\beta$,
$\gamma_{1}+i\gamma_{2}=\gamma$ and $u_{1}+iv_{2}=w$ \cite{comh}.
The set of equivalence transformations used to obtain the
semi-invariants of system (\ref{chypsysgen}) had a specific
Cauchy-Riemann (CR) structure. Therefore the deduced invariants also
satisfy the CR-equations. These quantities were also found to
correspond to their complex analogues given in \cite{ibrlp,lap},
associated with the base complex equation by means of a change of
the complex dependent variable. All the invariants of the scalar
linear hyperbolic equation are presented in \cite{fmjhs,ibrc}.

Semi-invariants of the linear parabolic equation via a
transformation of only the dependent variable have been derived
\cite{ibrlp}. Semi-invariants by a transformation of the independent
variables are given in \cite{mah,ibrb} and joint invariants are
found as well \cite{mah,ibrb,jf}. Laplace-type invariants of the
linear parabolic equation have been extended to Ibragimov-type
semi-invariants for a system of two parabolic-type PDEs with CSA
\cite{psysp}. Complex procedure imposes the CR-structure on the
invariants derived when applied on systems obtained by complex PDEs.
Though complex and real symmetry analysis have been employed to find
the semi-invariants of systems of parabolic and hyperbolic-type
PDEs, these results have not been compared earlier in order to bring
out the significance of CSA. By comparing the invariants associated
with the system of two hyperbolic-type PDEs with those deduced by
adopting real procedure, it is shown that CSA goes beyond the real
procedure and yields a few new invariant quantities different from
those obtained by real symmetry method developed for systems.

The plan of the paper is as follows. The second section is on
preliminaries where infinitesimal method is demonstrated. The
subsequent section contains derivation of the invariants of a system
of hyperbolic-type equations by real symmetry analysis. The fourth
section is on obtaining invariants for the same system by the
complex procedure and the comparison of these with the invariants
obtained in the third section. Concluding remarks are presented in
the last section.

\section{Preliminaries}
Semi-invariants of the linear hyperbolic equation
\begin{eqnarray}
w_{z_{1}z_{2}}+\alpha(z_{1},z_{2})w_{z_{1}}+\beta(z_{1},z_{2})w_{z_{2}}+\gamma(z_{1},z_{2})w=0,\label{hprblcsca}
\end{eqnarray}
under a transformation of (only) the dependent variables
\begin{eqnarray}
w(z_{1},z_{2})=\sigma(z_{1},z_{2})u(z_{1},z_{2}),\label{deptrans}
\end{eqnarray}
are derived in \cite{ibrlp,lap}. The transformations (\ref{deptrans})
correspond to the following infinitesimal change of the dependent
variables
\begin{eqnarray}
w=u+\epsilon\eta(z_{1},z_{2})u,\label{infdeptrans}
\end{eqnarray}
which leads to a generator of the form
\begin{eqnarray}
\textbf{Z}=\eta_{z_{2}}\partial_{\alpha}+\eta_{z_{1}}\partial_{\beta}+(\eta_{z_{1}z_{2}}+\alpha\eta_{z_{1}}+\beta\eta_{z_{2}})\partial_{\gamma},\label{depgen}
\end{eqnarray}
where $\eta_{z_{1}},~\eta_{z_{2}}$ denotes partial derivatives of
$\eta$, i.e. $\frac{\partial\eta}{\partial
z_1},~\frac{\partial\eta}{\partial z_2}$ and
$\partial_\alpha=\frac{\partial}{\partial\alpha}$, that is readable
from the transformed hyperbolic equation after implementing the
infinitesimal change of the dependent variable (\ref{infdeptrans})
on equation (\ref{hprblcsca}). The following first order
semi-invariants (called Laplace invariants) were reported in
\cite{ibrlp}
\begin{eqnarray}
h=\alpha_{z_{1}}+\alpha\beta-\gamma,~~k=\beta_{z_{2}}+\alpha\beta-\gamma,\label{hkinv}
\end{eqnarray}
by applying the infinitesimal method. The first extension of
(\ref{depgen}) acts on
$J(\alpha,\beta,\gamma,\alpha_{\rho},\beta_{\rho},\gamma_{\rho})$,
where $\rho=z_{1},z_{2}$, to reveal invariants (\ref{hkinv}).
Similarly, a change of the independent variables
\begin{eqnarray}
z_{1}=\phi(t),~~z_{2}=\psi(x),\label{indeptrans}
\end{eqnarray}
which can be written in the infinitesimal form as
\begin{eqnarray}
z_{1}=t+\epsilon\xi_{1}(t),~~~z_{2}=x+\epsilon\xi_{2}(x),\label{infindeptrans}
\end{eqnarray}
maps the linear hyperbolic equation (\ref{hprblcsca}) to an equation
of the same family, i.e. transformations (\ref{indeptrans}) preserve
the form, linearity and homogeneity of the hyperbolic PDE
(\ref{hprblcsca}). This change of (only) the independent variables
leads to a generator
\begin{eqnarray}
\textbf{Z}=\xi_{1}\partial_{t}+\xi_{2}\partial_{x}-\alpha\xi_{2,x}\partial_{\alpha}-\beta\xi_{1,t}\partial_{\beta}-\gamma(\xi_{1,t}+\xi_{2,x})\partial_{\gamma},\label{indepgen}
\end{eqnarray}
i.e. the operator obtained after transforming the linear hyperbolic
equation according to (\ref{infindeptrans}) and reading it from the
transformed new coefficients. Applying this generator on
$J(\alpha,\beta,\gamma)$ yields a zeroth order semi-invariant
\begin{eqnarray}
I_{1}=\frac{\gamma}{\alpha\beta}.
\end{eqnarray}
Further, in order to find the first order semi-invariants one needs
to engage the first extension of (\ref{indepgen}) on
$J(\alpha,\beta,\gamma,\alpha_{\varrho},\beta_{\varrho},\gamma_{\varrho})$,
where $\varrho=t,x$. This application results in a system of linear
PDEs that solves to the following invariant quantities
\begin{eqnarray}
I_{2}=\frac{\alpha\beta}{\alpha_{t}},~~I_{3}=\frac{\beta_{x}}{\alpha_{t}},~~I_{4}=\frac{\gamma}{\alpha_{t}},~~I_{5}=\frac{\alpha(\beta\gamma_{t}-\gamma\beta_{t})}{\beta\alpha^{2}_{t}},~~I_{6}=\frac{\alpha\gamma_{x}-\gamma\alpha_{x}}{\alpha^{2}\alpha_{t}}.
\end{eqnarray}
Further, the joint invariants of the hyperbolic equation have been
deduced in \cite{fmjhs,ibrc}. These invariants are derived by
applying an operator of the form (\ref{indepgen}) on the space of
the Laplace semi-invariants $h$ and $k$ given in (\ref{hkinv}).
Therefore, the first task is to transform the operator
(\ref{indepgen}) in the variables $h,~k$ and then apply it on
functions of these variables $J(h,k)$ and their derivatives
$J(h,k,h_{\varrho},k_{\varrho})$ and so on, in order to get the
zeroth, first and higher order \emph{joint invariants} of the linear
hyperbolic equation, respectively. After writing the generator
(\ref{indepgen}) in the space of the Laplace invariants $h,~k$, i.e.
\begin{eqnarray}
\textbf{Z}=\textbf{Z}(h)\partial_{h}+\textbf{Z}(k)\partial_{k},\label{tgensca}
\end{eqnarray}
its first extension reads as
\begin{eqnarray}
\textbf{Z}^{[1]}=\xi_{1}(t)\partial_{t}+\xi_{2}(x)\partial_{x}-(\xi_{1,t}+\xi_{2,x})h\partial_{h}-(\xi_{1,t}+\xi_{2,x})k\partial_{k}-(\xi_{1,tt}h+2\xi_{1,t}h_{t}+\xi_{2,x}h_{t})\partial_{h_{t}}\nonumber\\
-(\xi_{1,t}h_{x}+\xi_{2,xx}h+2\xi_{2,x}h_{x})\partial_{h_{x}}-(\xi_{1,tt}k+2\xi_{1,t}k_{t}+\xi_{2,x}k_{t})\partial_{k_{t}}-(\xi_{1,t}k_{x}+\xi_{2,xx}k~~~~~~~~\nonumber\\
+2\xi_{2,x}k_{x})\partial_{k_{x}}.~~~~~~~~~~~~~~~~~~~~~~~~~~~~~~~~~~~~~~~~~~~~~~~~~~~~~~~~~~~~~~~~~~~~~~~~~~~~~~~~~~~~~~~~~~~\label{trgenhksca}
\end{eqnarray}
It yields the following joint invariants \cite{fmjhs,ibrc}
\begin{eqnarray}
J_{1}=\frac{k}{h},~~~~~~~~~~~~~~~~~~~~~~~~~~~~~~~~~~~~~~~~~~~~~~~~~~~~~~~~~~~~~~\nonumber\\
J_{2}=\frac{(hk_{,t}-kh_{,t})(hk_{,x}-kh_{,x})}{h^5},~~~~~~~~~~~~~~~~~~~~~~~~~~~~~~~\nonumber\\
J_{3}=\frac{kh_{,tx}+hk_{,tx}-h_{,t}k_{,x}-h_{,x}k_{,t}}{h^3},~~~~~~~~~~~~~~~~~~~~~~~~~~~\nonumber\\
J_{4}=\frac{(hk_{,x}-kh_{,x})^{2}(hkh_{,tt}-h^2k_{,tt}-3kh_{,t}^{2}+3hh_{,t}k_{,t})}{h^9},~~\nonumber\\
J_{5}=\frac{(hk_{,t}-kh_{,t})^{2}(hkh_{,xx}-h^2k_{,xx}-3kh_{,x}^{2}+3hh_{,x}k_{,x})}{h^9},\nonumber\\
J_{6}=\frac{k(hh_{,tx}-h_{,t}h_{,x})}{h^4},~~~~~~~~~~~~~~~~~~~~~~~~~~~~~~~~~~~~~~~~~~~
\end{eqnarray}
of the scalar linear hyperbolic equation (\ref{hprblcsca}).

\section{Invariants of a system of two hyperbolic-type equations by real procedure}
Invariants of a system of two linear hyperbolic PDEs
(\ref{rehypsysgen}) have been determined by using the infinitesimal
method \cite{hsysp}. The derivation of these invariants starts with
determination of the equivalence transformations that map the system
of two linear hyperbolic equations into a system of the same form,
with different coefficients. The generators associated with these
infinitesimal transformations are then applied to obtain invariants
of (\ref{rehypsysgen}). Derivation of the invariants associated with
the subclass (\ref{chypsysgen}) of the system of hyperbolic
equations (\ref{rehypsysgen}) is presented in the remaining part of
this section, by real infinitesimal method.

The system of two hyperbolic-type PDEs (\ref{chypsysgen}) is
obtainable from a hyperbolic PDE with two independent variables,
when the dependent variable of equation (\ref{scacomhyp}) is
considered complex. Such systems have a CR-structure due to this
correspondence. Thus they are said to be \emph{CR-structured
systems}. The group of equivalence transformations associated with
(\ref{chypsysgen}) is obtained when the following generator
\begin{equation}
\begin{tabular}{l}
$\textbf{Z}=\xi^{1}\partial_{t}+\xi^{2}\partial_{x}+\eta^{1}\partial_{u}+\eta^{2}\partial_{v}+\eta^{1}_{t}\partial_{u_{t}}+\eta^{1}_{x}\partial_{u_{x}}+\eta^{2}_{t}\partial_{v_{t}}+\eta^{2}_{x}\partial_{v_{x}}+\eta^{1}_{tx}\partial_{u_{tx}}+\eta^{2}_{tx}\partial_{v_{tx}}$\\
$+\mu^{11}\partial_{\alpha_{1}}+\mu^{12}\partial_{\alpha_{2}}+\mu^{21}\partial_{\beta_{1}}+\mu^{22}\partial_{\beta_{2}}+\mu^{31}\partial_{\gamma_{1}}+\mu^{32}\partial_{\gamma_{2}},$\label{cgenhyp}
\end{tabular}
\end{equation}
acts on both the equations of the system (\ref{chypsysgen}). Where
$\xi^{\kappa}$, $\eta^{\kappa}$ are functions of $(t,x,u,v)$ and
$\mu^{1\kappa},\mu^{2\kappa}$, $\mu^{3\kappa}$, $\kappa=1,2,$ are
functions of
$(t,x,u,v,\alpha_{\kappa},\beta_{\kappa},\gamma_{\kappa})$, and
\begin{eqnarray}
&&\eta^{1}_{t}=D_{t}(\eta^{1})-u_{t}D_{t}(\xi^{1})-u_{x}D_{t}(\xi^{2}),\nonumber\\
&&\eta^{2}_{t}=D_{t}(\eta^{2})-v_{t}D_{t}(\xi^{1})-v_{x}D_{t}(\xi^{2}),\nonumber\\
&&\eta^{1}_{x}=D_{x}(\eta^{1})-u_{t}D_{x}(\xi^{1})-u_{x}D_{x}(\xi^{2}),\nonumber\\
&&\eta^{2}_{x}=D_{x}(\eta^{2})-v_{t}D_{x}(\xi^{1})-v_{x}D_{x}(\xi^{2}),
\end{eqnarray}
with $D_{t}$ and $D_{x}$ as total derivatives with respect to $t$
and $x$. Implication of (\ref{cgenhyp}) yields a system of linear
PDEs that leads to
\begin{equation}
\begin{tabular}{l}
$\xi_{1}=F_{1}(t),~~\xi_{2}=F_{2}(x)$,\\
$\mu^{11}=-F_{3,x}-\alpha_{1}F_{2,x},~~\mu^{12}=F_{4,x}-\alpha_{2}F_{2,x}$,\\
$\mu^{21}=-F_{3,t}-\beta_{1}F_{1,t},~~\mu^{22}=F_{4,t}-\beta_{2}F_{1,t}$,\\
$\mu^{31}=-F_{3,tx}-\alpha_{1}F_{3,t}-\alpha_{2}F_{4,t}-\beta_{1}F_{3,x}-\beta_{2}F_{4,x}-\gamma_{1}(F_{1,t}+F_{2,x})$,\\
$\mu^{32}=F_{4,tx}+\alpha_{1}F_{4,t}-\alpha_{2}F_{3,t}+\beta_{1}F_{4,x}-\beta_{2}F_{3,x}-\gamma_{2}(F_{1,t}+F_{2,x})$,\label{genrtreal}
\end{tabular}
\end{equation}
where $F_{3}$ and $F_{4}$ depends on $(t,x)$.

The first order semi-invariants
\begin{eqnarray}
h^{^{r}}_{1}=\alpha_{1,t}+\alpha_{1}\beta_{1}-\alpha_{2}\beta_{2}-\gamma_{1},\nonumber\\
h^{^{r}}_{2}=\alpha_{2,t}+\alpha_{1}\beta_{2}+\alpha_{2}\beta_{1}-\gamma_{2},\nonumber\\
k^{^{r}}_{1}=\beta_{1,x}+\alpha_{1}\beta_{1}-\alpha_{2}\beta_{2}-\gamma_{1},\nonumber\\
k^{^{r}}_{2}=\beta_{2,x}+\alpha_{1}\beta_{2}+\alpha_{2}\beta_{1}-\gamma_{2},\label{rhksys}
\end{eqnarray}
associated with the system (\ref{chypsysgen}) due to a change of
(only) dependent variables
\begin{eqnarray}
\eta^{1}=F_{3}u+F_{4}v,~~\eta^{2}=F_{3}v-F_{4}u,
\end{eqnarray}
are obtained by employing the generator
\begin{eqnarray}
&&\textbf{X}=-F_{3,x}\partial_{\alpha_{1}}+F_{4,x}\partial_{\alpha_{2}}-F_{3,t}\partial_{\beta_{1}}+F_{4,t}\partial_{\beta_{2}}-(F_{3,tx}+\alpha_{1}F_{3,t}+\alpha_{2}F_{4,t}\nonumber\\
&&~~~~+\beta_{1}F_{3,x}+\beta_{2}F_{4,x})\partial_{\gamma_{1}}+(F_{4,tx}+\alpha_{1}F_{4,t}-\alpha_{2}F_{3,t}+\beta_{1}F_{4,x}-\beta_{2}F_{3,x})\partial_{\gamma_{2}}.\label{cgenhypd}
\end{eqnarray}
It is extracted from (\ref{genrtreal}) by considering
$F_{1}(t)=F_{2}(x)=0$, and $F_{3}(t,x),~F_{4}(t,x)$ as arbitrary
functions of their arguments. When it acts on
$J(\alpha_{\kappa},\beta_{\kappa},\gamma_{\kappa},\alpha_{\kappa,t},\beta_{\kappa,t},\gamma_{\kappa,t},\alpha_{\kappa,x},\beta_{\kappa,x},\gamma_{\kappa,x})$
we obtain (\ref{rhksys}) by solving the resulting linear system of
PDEs.

Considering only a change of the independent variables, i.e.,
keeping $F_{1}(t),~F_{2}(x)$ as arbitrary functions of their
arguments and $F_{3}(t,x)=F_{4}(t,x)=0$, leads to an infinitesimal
generator
\begin{eqnarray}
&&\textbf{Z}_{I}=F_{1}(t)\partial_{t}+F_{2}(x)\partial_{x}-\alpha_{1}F_{2,x}\partial_{\alpha_{1}}-\alpha_{2}F_{2,x}\partial_{\alpha_{2}}-\beta_{1}F_{1,t}\partial_{\beta_{1}}-\beta_{2}F_{1,t}\partial_{\beta_{2}}\nonumber\\
&&~~~~-\gamma_{1}(F_{1,t}+F_{2,x})\partial_{\gamma_{1}}-\gamma_{2}(F_{1,t}+F_{2,x})\partial_{\gamma_{2}}.\label{cgenhypI}
\end{eqnarray}
Applying it on $J(\alpha_{\kappa},\beta_{\kappa},\gamma_{\kappa})$
yields the following zeroth order invariants
\begin{eqnarray}
I^{^{r}}_{1}=\frac{\alpha_{2}}{\alpha_{1}},~~I^{^{r}}_{2}=\frac{\beta_{2}}{\beta_{1}},~~I^{^{r}}_{3}=\frac{\gamma_{1}}{\alpha_{1}\beta_{1}},~~I^{^{r}}_{4}=\frac{\gamma_{2}}{\alpha_{1}\beta_{1}}.\label{0orderindep}
\end{eqnarray}
Further, the first order invariants are obtained when the once
extended generator (\ref{cgenhypI}) acts on
$J(\alpha_{\kappa},\beta_{\kappa},\gamma_{\kappa},\alpha_{\kappa,t},\beta_{\kappa,t},\gamma_{\kappa,t},\alpha_{\kappa,x},\beta_{\kappa,x},\gamma_{\kappa,x})$,
gives the following quantities
\begin{eqnarray}
I^{^{r}}_{5}=\frac{\alpha_{1,t}}{\alpha_{1}\beta_{1}},~~I^{^{r}}_{6}=\frac{\alpha_{2,t}}{\alpha_{1}\beta_{1}},~~I^{^{r}}_{7}=\frac{\beta_{1,x}}{\alpha_{1}\beta_{1}},~~I^{^{r}}_{8}=\frac{\beta_{2,x}}{\alpha_{1}\beta_{1}},~~~~~~~~~~~~~~~~~~~~~~~~~\nonumber\\
I^{^{r}}_{9}=\frac{\beta_{1}\beta_{2,t}-\beta_{2}\beta_{1,t}}{b^{3}_{1}},~~I^{^{r}}_{10}=\frac{\beta_{1}\gamma_{1,t}-\gamma_{1}\beta_{1,t}}{\alpha_{1}\beta^{3}_{1}},~~I^{^{r}}_{11}=\frac{\beta_{1}\gamma_{2,t}-\gamma_{2}\beta_{1,t}}{\alpha_{1}\beta^{3}_{1}},~~~~~\nonumber\\
I^{^{r}}_{12}=\frac{\alpha_{1}\alpha_{2,x}-\alpha_{2}\alpha_{1,x}}{\alpha^{3}_{1}},~~I^{^{r}}_{13}=\frac{\alpha_{1}\gamma_{1,x}-\gamma_{1}\alpha_{1,x}}{\alpha^{3}_{1}\beta_{1}},~~I^{^{r}}_{14}=\frac{\alpha_{1}\gamma_{2,x}-\gamma_{2}\alpha_{1,x}}{\alpha^{3}_{1}\beta_{1}},
\end{eqnarray}
including the four zeroth order invariants (\ref{0orderindep}).

The joint invariants of the system (\ref{chypsysgen})
\begin{eqnarray}
J^{^{r}}_{1}=\frac{h^{^{r}}_{2}}{h^{^{r}}_{1}},~~J^{^{r}}_{2}=\frac{k^{^{r}}_{1}}{h^{^{r}}_{1}},~~J^{^{r}}_{3}=\frac{k^{^{r}}_{2}}{h^{^{r}}_{1}},
\end{eqnarray}
are found when the following PDE
\begin{eqnarray}
h^{^{r}}_{1}\partial_{h^{^{r}}_{1}}+h^{^{r}}_{2}\partial_{h^{^{r}}_{2}}+k^{^{r}}_{1}\partial_{k^{^{r}}_{1}}+k^{^{r}}_{2}\partial_{k^{^{r}}_{2}}=0,
\end{eqnarray}
is solved. This equation appears due to action of the
infinitesimal generator (\ref{cgenhypI}) that is associated with the
change of the independent variables to the space of invariants
$h^{^{r}}_{\kappa},~k^{^{r}}_{\kappa}$.

\section{Invariants of a system of two hyperbolic-type equations by complex procedure}
Semi-invariants associated with a system of two hyperbolic-type
equations (\ref{chypsysgen}) that is obtained from a scalar linear
hyperbolic equation (\ref{scacomhyp}), are derived in this section
by complex methods. A few of the invariants presented here have
already been provided \cite{comh}, we demonstrate the complete
complex procedure involved to re-derive them. The generator of the
form (\ref{depgen}) associated with the equation (\ref{hprblcsca})
becomes complex due to the presence of the complex dependent
variable and the complex coefficients split (\ref{depgen}) into two
operators
\begin{eqnarray}
\textbf{X}_{1}=\eta_{1,z_{2}}\partial_{\alpha_{1}}+\eta_{2,z_{2}}\partial_{\alpha_{2}}+\eta_{1,z_{1}}\partial_{\beta_{1}}+\eta_{2,z_{1}}\partial_{\beta_{2}}+(\eta_{1,z_{1}z_{2}}+\alpha_{1}\eta_{1,z_{1}}-\alpha_{2}\eta_{2,z_{1}}+\beta_{1}\eta_{1,z_{2}}\nonumber\\
-\beta_{2}\eta_{2,z_{2}})\partial_{\gamma_{1}}+(\eta_{2,z_{1}z_{2}}+\alpha_{2}\eta_{1,z_{1}}+\alpha_{1}\eta_{2,z_{1}}+\beta_{2}\eta_{1,z_{2}}+\beta_{1}\eta_{2,z_{2}})\partial_{\gamma_{2}},~~~~~~~~~~~~~~~~~\label{depop1}
\end{eqnarray}
\begin{eqnarray}
\textbf{X}_{2}=\eta_{2,z_{2}}\partial_{\alpha_{1}}-\eta_{1,z_{2}}\partial_{\alpha_{2}}+\eta_{2,z_{1}}\partial_{\beta_{1}}-\eta_{1,z_{1}}\partial_{\beta_{2}}+(\eta_{2,z_{1}z_{2}}+\alpha_{2}\eta_{1,z_{1}}+\alpha_{1}\eta_{2,z_{1}}+\beta_{2}\eta_{1,z_{2}}\nonumber\\
+\beta_{1}\eta_{2,z_{2}})\partial_{\gamma_{1}}-(\eta_{1,z_{1}z_{2}}+\alpha_{1}\eta_{1,z_{1}}-\alpha_{2}\eta_{2,z_{1}}+\beta_{1}\eta_{1,z_{2}}-\beta_{2}\eta_{2,z_{2}})\partial_{\gamma_{2}}.~~~~~~~~~~~~~~~~~\label{depop2}
\end{eqnarray}
There are four first order semi-invariants
\begin{eqnarray}
h_{1}=\alpha_{1,z_{1}}+\alpha_{1}\beta_{1}-\alpha_{2}\beta_{2}-\gamma_{1},\nonumber\\
h_{2}=\alpha_{2,z_{1}}+\alpha_{2}\beta_{1}+\alpha_{1}\beta_{2}-\gamma_{2},\nonumber\\
k_{1}=\beta_{1,z_{2}}+\alpha_{1}\beta_{1}-\alpha_{2}\beta_{2}-\gamma_{1},\nonumber\\
k_{2}=\beta_{2,z_{2}}+\alpha_{2}\beta_{1}+\alpha_{1}\beta_{2}-\gamma_{2},\label{fourinv}
\end{eqnarray}
that are found to be associated with the system (\ref{chypsysgen})
on employing the pair of operators (\ref{depop1}) and
(\ref{depop2}). These are exactly the same as represented by
$h^{^{r}}_{\kappa},~k^{^{r}}_{\kappa}$ in (\ref{rhksys}). Therefore,
in this case the real and complex procedures lead to the same
semi-invariants of the system (\ref{chypsysgen}). Notice that all
the four semi-invariants (\ref{fourinv}) are readable from the the
first order semi-invariants associated with the complex hyperbolic
linear equation (\ref{hprblcsca}) and satisfy
\begin{equation}
\begin{tabular}{l}
$\textbf{X}^{[1]}_{1}h_{1}\mid_{_{_{_{h_{1}=0}}}}=\textbf{X}^{[1]}_{2}h_{2}\mid_{_{_{_{h_{2}=0}}}}=\textbf{X}^{[1]}_{1}k_{1}\mid_{_{_{_{k_{1}=0}}}}=\textbf{X}^{[1]}_{2}k_{2}\mid_{_{_{_{k_{2}=0}}}}=0$.
\end{tabular}
\end{equation}
The linear combination $\textbf{X}_{3}$ of both the operators $\textbf{X}_{1}$ and $\textbf{X}_{2}$ results in
the following relations
\begin{equation}
\begin{tabular}{l}
$\textbf{X}^{[1]}_{3}h_{1}\mid_{_{_{_{h_{1}=0}}}}=\textbf{X}^{[1]}_{3}h_{2}\mid_{_{_{_{h_{2}=0}}}}=\textbf{X}^{[1]}_{3}k_{1}\mid_{_{_{_{k_{1}=0}}}}=\textbf{X}^{[1]}_{3}k_{2}\mid_{_{_{_{k_{2}=0}}}}=0$.
\end{tabular}
\end{equation}
The semi-invariants of the system of two hyperbolic-type PDEs under
a transformation of the independent variables are
\begin{eqnarray*}
I^{^{c}}_{1}=\frac{(\alpha_{1}\gamma_{1}+\alpha_{2}\gamma_{2})\beta_{1}+(\alpha_{1}\gamma_{2}-\alpha_{2}\gamma_{1})\beta_{2}}{(\alpha^{2}_{1}+\alpha^{2}_{2})(\beta^{2}_{1}+\beta^{2}_{2})},\nonumber\\
I^{^{c}}_{2}=\frac{(\alpha_{1}\gamma_{2}-\alpha_{2}\gamma_{1})\beta_{1}-(\alpha_{1}\gamma_{1}+\alpha_{2}\gamma_{2})\beta_{2}}{(\alpha^{2}_{1}+\alpha^{2}_{2})(\beta^{2}_{1}+\beta^{2}_{2})},
\end{eqnarray*}
\begin{eqnarray*}
I^{^{c}}_{3}=\frac{(\alpha_{1}\beta_{1}-\alpha_{2}\beta_{2})\alpha_{1,t}+(\alpha_{2}\beta_{1}+\alpha_{1}\beta_{2})\alpha_{2,t}}{\alpha^{2}_{1,t}+\alpha^{2}_{2,t}},\nonumber\\
I^{^{c}}_{4}=\frac{(\alpha_{2}\beta_{1}+\alpha_{1}\beta_{2})\alpha_{1,t}-(\alpha_{1}\beta_{1}-\alpha_{2}\beta_{2})\alpha_{2,t}}{\alpha^{2}_{1,t}+\alpha^{2}_{2,t}},
\end{eqnarray*}
\begin{eqnarray*}
I^{^{c}}_{5}=\frac{\alpha_{1,t}\beta_{1,x}+\alpha_{2,t}\beta_{2,x}}{\alpha^{2}_{1,t}+\alpha^{2}_{2,t}},~~~I^{^{c}}_{6}=\frac{\alpha_{1,t}\beta_{2,x}-\alpha_{2,t}\beta_{1,x}}{\alpha^{2}_{1,t}+\alpha^{2}_{2,t}},
\end{eqnarray*}
\begin{eqnarray*}
I^{^{c}}_{7}=\frac{\alpha_{1,t}\gamma_{1}+\alpha_{2,t}\gamma_{2}}{\alpha^{2}_{1,t}+\alpha^{2}_{2,t}},~~~I^{^{c}}_{8}=\frac{\alpha_{1,t}\gamma_{2}-\alpha_{2,t}\gamma_{1}}{\alpha^{2}_{1,t}+\alpha^{2}_{2,t}},
\end{eqnarray*}
\begin{equation}
\begin{tabular}{l}
$I^{^{c}}_{9}=\frac{(\alpha^{2}_{1,t}-\alpha^{2}_{2,t})}{(\alpha^{2}_{1,t}+\alpha^{2}_{2,t})^{2}(\beta_{1}^{2}+\beta_{2}^{2})}[\alpha_{1}\beta_{1}(\beta_{1}\gamma_{1,t}-\beta_{2}\gamma_{2,t}-\gamma_{1}\beta_{1,t}+\gamma_{2}\beta_{2,t})-\alpha_{2}\beta_{1}(\beta_{2}\gamma_{1,t}+\beta_{1}\gamma_{2,t}$\\
$-\gamma_{2}\beta_{1,t}-\gamma_{1}\beta_{2,t})+\alpha_{2}\beta_{2}(\beta_{1}\gamma_{1,t}-\beta_{2}\gamma_{2,t}-\gamma_{1}\beta_{1,t}+\gamma_{2}\beta_{2,t})+\alpha_{1}\beta_{2}(\beta_{2}\gamma_{1,t}+\beta_{1}\gamma_{2,t}$\\
$-\gamma_{2}\beta_{1,t}-\gamma_{1}\beta_{2,t})]+\frac{2\alpha_{1,t}\alpha_{2,t}}{(\alpha^{2}_{1,t}+\alpha^{2}_{2,t})^{2}(\beta_{1}^{2}+\beta_{2}^{2})}[\alpha_{2}\beta_{1}(\beta_{1}\gamma_{1,t}-\beta_{2}\gamma_{2,t}-\gamma_{1}\beta_{1,t}+\gamma_{2}\beta_{2,t})$\\
$+\alpha_{1}\beta_{1}(\beta_{2}\gamma_{1,t}+\beta_{1}\gamma_{2,t}-\gamma_{2}\beta_{1,t}-\gamma_{1}\beta_{2,t})-\alpha_{1}\beta_{2}(\beta_{1}\gamma_{1,t}-\beta_{2}\gamma_{2,t}-\gamma_{1}\beta_{1,t}+\gamma_{2}\beta_{2,t})$\\
$+\alpha_{2}\beta_{2}(\beta_{2}\gamma_{1,t}+\beta_{1}\gamma_{2,t}-\gamma_{2}\beta_{1,t}-\gamma_{1}\beta_{2,t})],$
\end{tabular}
\end{equation}
\begin{equation}
\begin{tabular}{l}
$I^{^{c}}_{10}=\frac{(\alpha^{2}_{1,t}-\alpha^{2}_{2,t})}{(\alpha^{2}_{1,t}+\alpha^{2}_{2,t})^{2}(\beta_{1}^{2}+\beta_{2}^{2})}[\alpha_{2}\beta_{1}(\beta_{1}\gamma_{1,t}-\beta_{2}\gamma_{2,t}-\gamma_{1}\beta_{1,t}+\gamma_{2}\beta_{2,t})+\alpha_{1}\beta_{1}(\beta_{2}\gamma_{1,t}+\beta_{1}\gamma_{2,t}$\\
$-\gamma_{2}\beta_{1,t}-\gamma_{1}\beta_{2,t})-\alpha_{1}\beta_{2}(\beta_{1}\gamma_{1,t}-\beta_{2}\gamma_{2,t}-\gamma_{1}\beta_{1,t}+\gamma_{2}\beta_{2,t})+\alpha_{2}\beta_{2}(\beta_{2}\gamma_{1,t}+\beta_{1}\gamma_{2,t}$\\
$-\gamma_{2}\beta_{1,t}-\gamma_{1}\beta_{2,t})]-\frac{2\alpha_{1,t}\alpha_{2,t}}{(\alpha^{2}_{1,t}+\alpha^{2}_{2,t})^{2}(\beta_{1}^{2}+\beta_{2}^{2})}[\alpha_{1}\beta_{1}(\beta_{1}\gamma_{1,t}-\beta_{2}\gamma_{2,t}-\gamma_{1}\beta_{1,t}+\gamma_{2}\beta_{2,t})$\\
$-\alpha_{2}\beta_{1}(\beta_{2}\gamma_{1,t}+\beta_{1}\gamma_{2,t}-\gamma_{2}\beta_{1,t}-\gamma_{1}\beta_{2,t})+\alpha_{2}\beta_{2}(\beta_{1}\gamma_{1,t}-\beta_{2}\gamma_{2,t}-\gamma_{1}\beta_{1,t}+\gamma_{2}\beta_{2,t})$\\
$+\alpha_{1}\beta_{2}(\beta_{2}\gamma_{1,t}+\beta_{1}\gamma_{2,t}-\gamma_{2}\beta_{1,t}-\gamma_{1}\beta_{2,t})],$
\end{tabular}
\end{equation}
\begin{equation}
\begin{tabular}{l}
$I^{^{c}}_{11}=\frac{(\alpha^{2}_{1}-\alpha^{2}_{2})}{(\alpha^{2}_{1,t}+\alpha^{2}_{2,t})(\alpha_{1}^{2}+\alpha_{2}^{2})^{2}}[(\alpha_{1}\gamma_{1,x}-\alpha_{2}\gamma_{2,x}-\gamma_{1}\alpha_{1,x}+\gamma_{2}\alpha_{2,x})\alpha_{1,t}+(\alpha_{2}\gamma_{1,x}+\alpha_{1}\gamma_{2,x}$\\
$-\gamma_{2}\alpha_{1,x}-\gamma_{1}\alpha_{2,x})\alpha_{2,t}]+\frac{2\alpha_{1}\alpha_{2}}{(\alpha^{2}_{1,t}+\alpha^{2}_{2,t})(\alpha_{1}^{2}+\alpha_{2}^{2})^{2}}[(\alpha_{2}\gamma_{1,x}+\alpha_{1}\gamma_{2,x}-\gamma_{2}\alpha_{1,x}-\gamma_{1}\alpha_{2,x})\alpha_{1,t}$\\
$-(\alpha_{1}\gamma_{1,x}-\alpha_{1}\gamma_{2,x}-\gamma_{1}\alpha_{1,x}+\gamma_{2}\alpha_{2,x})\alpha_{2,t}],$
\end{tabular}
\end{equation}
\begin{equation}
\begin{tabular}{l}
$I^{^{c}}_{12}=\frac{(\alpha^{2}_{1}-\alpha^{2}_{2})}{(\alpha^{2}_{1,t}+\alpha^{2}_{2,t})(\alpha_{1}^{2}+\alpha_{2}^{2})^{2}}[(\alpha_{2}\gamma_{1,x}+\alpha_{1}\gamma_{2,x}-\gamma_{2}\alpha_{1,x}-\gamma_{1}\alpha_{2,x})\alpha_{1,t}-(\alpha_{1}\gamma_{1,x}-\alpha_{2}\gamma_{2,x}$\\
$-\gamma_{1}\alpha_{1,x}+\gamma_{2}\alpha_{2,x})\alpha_{2,t}]-\frac{2\alpha_{1}\alpha_{2}}{(\alpha^{2}_{1,t}+\alpha^{2}_{2,t})(\alpha_{1}^{2}+\alpha_{2}^{2})^{2}}[(\alpha_{1}\gamma_{1,x}-\alpha_{2}\gamma_{2,x}-\gamma_{1}\alpha_{1,x}+\gamma_{2}\alpha_{2,x})\alpha_{1,t}$\\
$+(\alpha_{2}\gamma_{1,x}+\alpha_{1}\gamma_{2,x}-\gamma_{2}\alpha_{1,x}-\gamma_{1}\alpha_{2,x})\alpha_{2,t}].$
\end{tabular}
\end{equation}
The correspondence of these semi-invariants of independent variables with the system of the
hyperbolic-type equations is established due to the following
operators
\begin{eqnarray}
\textbf{X}_{1}=2\xi_{1}\partial_{t}+2\xi_{2}\partial_{x}-\alpha_{1}\xi_{2,x}\partial_{\alpha_{1}}-\alpha_{2}\xi_{2,x}\partial_{\alpha_{2}}-\beta_{1}\xi_{1,t}\partial_{\beta_{1}}-\beta_{2}\xi_{1,t}\partial_{\beta_{2}}-\gamma_{1}(\xi_{1,t}\nonumber\\
+\xi_{2,x})\partial_{\gamma_{1}}-\gamma_{2}(\xi_{1,t}+\xi_{2,x})\partial_{\gamma_{2}},~~~~~~~~~~~~~~~~~~~~~~~~~~~~~~~~~~~~~~~~~~~~~~~~~~~~\label{indepop1}
\end{eqnarray}
\begin{eqnarray}
\textbf{X}_{2}=-\alpha_{2}\xi_{2,x}\partial_{\alpha_{1}}+\alpha_{1}\xi_{2,x}\partial_{\alpha_{2}}-\beta_{2}\xi_{1,t}\partial_{\beta_{1}}+\beta_{1}\xi_{1,t}\partial_{\beta_{2}}-\gamma_{2}(\xi_{1,t}+\xi_{2,x})\partial_{\gamma_{1}}~~~~~\nonumber\\
+\gamma_{1}(\xi_{1,t}+\xi_{2,x})\partial_{\gamma_{2}},~~~~~~~~~~~~~~~~~~~~~~~~~~~~~~~~~~~~~~~~~~~~~~~~~~~~~~~~~~~~~~~\label{indepop2}
\end{eqnarray}
that are the real and imaginary parts of the complex generator
(\ref{indepgen}). Using these operators it is observed that
\begin{equation}
\begin{tabular}{l}
$\textbf{X}^{[1]}_{1}I^{^{c}}_{1}\mid_{_{_{_{I^{^{c}}_{1}=0}}}}=\textbf{X}^{[1]}_{2}I^{^{c}}_{2}\mid_{_{_{_{I^{^{c}}_{2}=0}}}}=\textbf{X}^{[1]}_{1}I^{^{c}}_{3}\mid_{_{_{_{I^{^{c}}_{3}=0}}}}=\textbf{X}^{[1]}_{2}I^{^{c}}_{4}\mid_{_{_{_{I^{^{c}}_{3}=I^{^{c}}_{4}=0}}}}=0$,\\
$\textbf{X}^{[1]}_{1}I^{^{c}}_{5}\mid_{_{_{_{I^{^{c}}_{5}=0}}}}=\textbf{X}^{[1]}_{2}I^{^{c}}_{6}\mid_{_{_{_{I^{^{c}}_{5}=I^{^{c}}_{6}=0}}}}=\textbf{X}^{[1]}_{1}I^{^{c}}_{7}\mid_{_{_{_{I^{^{c}}_{7}=0}}}}=\textbf{X}^{[1]}_{2}I^{^{c}}_{8}\mid_{_{_{_{I^{^{c}}_{7}=I^{^{c}}_{8}=0}}}}=0$,\\
$\textbf{X}^{[1]}_{1}I^{^{c}}_{9}\mid_{_{_{_{I^{^{c}}_{9}=0}}}}=\textbf{X}^{[1]}_{2}I^{^{c}}_{10}\mid_{_{_{_{I^{^{c}}_{9}=I^{^{c}}_{10}=0}}}}=\textbf{X}^{[1]}_{1}I^{^{c}}_{11}\mid_{_{_{_{I^{^{c}}_{11}=0}}}}=\textbf{X}^{[1]}_{2}I^{^{c}}_{12}\mid_{_{_{_{I^{^{c}}_{11}=I^{^{c}}_{12}=0}}}}=0$.
\end{tabular}
\end{equation}
It is seen that the above invariants are complex splits of their real analogues (13).
Similarly, the linear combination of both $\textbf{X}_{1}$ and
$\textbf{X}_{2}$, if denoted by $\textbf{X}_{3}$, satisfy the
 relations
\begin{equation}
\begin{tabular}{l}
$\textbf{X}^{[1]}_{3}I^{^{c}}_{1}\mid_{_{_{_{I^{^{c}}_{1}=0}}}}=\textbf{X}^{[1]}_{3}I^{^{c}}_{2}\mid_{_{_{_{I^{^{c}}_{2}=0}}}}=\textbf{X}^{[1]}_{3}I^{^{c}}_{3}\mid_{_{_{_{I^{^{c}}_{3}=0}}}}=\textbf{X}^{[1]}_{3}I^{^{c}}_{4}\mid_{_{_{_{I^{^{c}}_{3}=I^{^{c}}_{4}=0}}}}=0$,\\
$\textbf{X}^{[1]}_{3}I^{^{c}}_{5}\mid_{_{_{_{I^{^{c}}_{5}=0}}}}=\textbf{X}^{[1]}_{3}I^{^{c}}_{6}\mid_{_{_{_{I^{^{c}}_{5}=I^{^{c}}_{6}=0}}}}=\textbf{X}^{[1]}_{3}I^{^{c}}_{7}\mid_{_{_{_{I^{^{c}}_{7}=I^{^{c}}_{8}=0}}}}=\textbf{X}^{[1]}_{3}I^{^{c}}_{8}\mid_{_{_{_{I^{^{c}}_{7}=I^{^{c}}_{8}=0}}}}=0$,\\
$\textbf{X}^{[1]}_{3}I^{^{c}}_{9}\mid_{_{_{_{I^{^{c}}_{9}=I^{^{c}}_{10}=0}}}}=\textbf{X}^{[1]}_{3}I^{^{c}}_{10}\mid_{_{_{_{I^{^{c}}_{9}=I^{^{c}}_{10}=0}}}}=\textbf{X}^{[1]}_{3}I^{^{c}}_{11}\mid_{_{_{_{I^{^{c}}_{11}=I^{^{c}}_{12}=0}}}}=\textbf{X}^{[1]}_{3}I^{^{c}}_{12}\mid_{_{_{_{I^{^{c}}_{11}=I^{^{c}}_{12}=0}}}}=0$.
\end{tabular}
\end{equation}
To work out the joint invariants of the coupled system of two
hyperbolic-type equations (\ref{chypsysgen}), the operators
(\ref{indepop1}) and (\ref{indepop2}) need to be transformed to the
space of invariants $h_{\kappa},~k_{\kappa}$. The same procedure was
adopted in \cite{fmjhs} before using the generator (\ref{indepgen})
in determining the joint invariants of the scalar linear hyperbolic
equation. The complex generator was transformed to $h$ and $k$,
i.e. to the space of the semi-invariants associated with the hyperbolic equation
under a change of the dependent variables. The procedure to
transform (\ref{indepop1}) and (\ref{indepop2}) to
$(h_{\kappa},k_{\kappa})-space$ starts with splitting
(\ref{tgensca}) when $\textbf{Z}(h)$ and $\textbf{Z}(k)$ are
taken as complex, i.e.
$\textbf{Z}(h)=\textbf{Z}(h)_{1}+i\textbf{Z}(h)_{2}$ and
$\textbf{Z}(k)=\textbf{Z}(k)_{1}+i\textbf{Z}(k)_{2}$. The real and
imaginary parts of (\ref{tgensca}) are
\begin{eqnarray}
\textbf{X}_{1}=\frac{1}{2}[\textbf{Z}(h)_{1}\partial_{h_{1}}+\textbf{Z}(h)_{2}\partial_{h_{2}}+\textbf{Z}(k)_{1}\partial_{k_{1}}+\textbf{Z}(k)_{2}\partial_{k_{2}}],\nonumber\\
\textbf{X}_{2}=\frac{1}{2}[\textbf{Z}(h)_{2}\partial_{h_{1}}-\textbf{Z}(h)_{1}\partial_{h_{2}}+\textbf{Z}(k)_{2}\partial_{k_{1}}-\textbf{Z}(k)_{1}\partial_{k_{2}}],\label{genops1}
\end{eqnarray}
where
\begin{eqnarray}
\textbf{Z}(h)_{1}=\textbf{X}_{1}h_{1}-\textbf{X}_{2}h_{2}=-(\xi_{1,t}+\xi_{2,x})h_{1},\nonumber\\
\textbf{Z}(h)_{2}=\textbf{X}_{2}h_{1}+\textbf{X}_{1}h_{2}=-(\xi_{1,t}+\xi_{2,x})h_{2},\nonumber\\
\textbf{Z}(k)_{1}=\textbf{X}_{1}k_{1}-\textbf{X}_{2}k_{2}=-(\xi_{1,t}+\xi_{2,x})k_{1},\nonumber\\
\textbf{Z}(k)_{2}=\textbf{X}_{2}k_{1}+\textbf{X}_{1}k_{2}=-(\xi_{1,t}+\xi_{2,x})k_{2}.\label{genops2}
\end{eqnarray}
Using (\ref{genops2}) in (\ref{genops1}) yields the following two
operators
\begin{eqnarray}
\textbf{X}_{1}=-\frac{(\xi_{1,t}+\xi_{2,x})}{2}[h_{1}\partial_{h_{1}}+h_{2}\partial_{h_{2}}+k_{1}\partial_{k_{1}}+k_{2}\partial_{k_{2}}],\nonumber\\
\textbf{X}_{2}=-\frac{(\xi_{1,t}+\xi_{2,x})}{2}[h_{2}\partial_{h_{1}}-h_{1}\partial_{h_{2}}+k_{2}\partial_{k_{1}}-k_{1}\partial_{k_{2}}],
\end{eqnarray}
that are the real and imaginary parts of the complex generator
(\ref{trgenhksca}). These operators are used to arrive at the joint
invariants for the system of two linear hyperbolic-type equations
(\ref{chypsysgen}). We have the following joint invariants
\begin{eqnarray*}
J_{11}=\frac{h_{1}k_{1}+h_{2}k_{2}}{k^{2}_{1}+k^{2}_{2}},\nonumber\\
J_{12}=\frac{h_{2}k_{1}-h_{1}k_{2}}{k^{2}_{1}+k^{2}_{2}},
\end{eqnarray*}
\begin{eqnarray*}
J_{13}=\frac{(h_{1}^{5}-10h_{1}^{3}h_{2}^{2}+5h_{1}h_{2}^{4})(h_{1}k_{1,t}-h_{2}k_{2,t}-k_{1}h_{1,t}+k_{2}h_{2,t})}{(h_{1}^{5}-10h_{1}^{3}h_{2}^{2}+5h_{1}h_{2}^{4})^{2}+(5h_{1}^{4}h_{2}-10h_{1}^{2}h_{2}^{3}+h_{2}^{5})^{2}}(h_{1}k_{1,x}-h_{2}k_{2,x}-k_{1}h_{1,x}+k_{2}h_{2,x})\nonumber\\
+\frac{(5h_{1}^{4}h_{2}-10h_{1}^{2}h_{2}^{3}+h_{2}^{5})(h_{2}k_{1,t}+h_{1}k_{2,t}-k_{2}h_{1,t}-k_{1}h_{2,t})}{(h_{1}^{5}-10h_{1}^{3}h_{2}^{2}+5h_{1}h_{2}^{4})^{2}+(5h_{1}^{4}h_{2}-10h_{1}^{2}h_{2}^{3}+h_{2}^{5})^{2}}(h_{1}k_{1,x}-h_{2}k_{2,x}-k_{1}h_{1,x}+k_{2}h_{2,x})\nonumber\\
+\frac{(5h_{1}^{4}h_{2}-10h_{1}^{2}h_{2}^{3}+h_{2}^{5})(h_{1}k_{1,t}-h_{2}k_{2,t}-k_{1}h_{1,t}+k_{2}h_{2,t})}{(h_{1}^{5}-10h_{1}^{3}h_{2}^{2}+5h_{1}h_{2}^{4})^{2}+(5h_{1}^{4}h_{2}-10h_{1}^{2}h_{2}^{3}+h_{2}^{5})^{2}}(h_{2}k_{1,x}+h_{1}k_{2,x}-k_{2}h_{1,x}-k_{1}h_{2,x})\nonumber\\
+\frac{(h_{1}^{5}-10h_{1}^{3}h_{2}^{2}+5h_{1}h_{2}^{4})(h_{2}k_{1,t}+h_{1}k_{2,t}-k_{1}h_{1,t}-k_{1}h_{2,t})}{(h_{1}^{5}-10h_{1}^{3}h_{2}^{2}+5h_{1}h_{2}^{4})^{2}+(5h_{1}^{4}h_{2}-10h_{1}^{2}h_{2}^{3}+h_{2}^{5})^{2}}(h_{2}k_{1,x}+h_{1}k_{2,x}-k_{2}h_{1,x}-k_{1}h_{2,x}),
\end{eqnarray*}
\begin{eqnarray*}
J_{14}=\frac{-(5h_{1}^{4}h_{2}-10h_{1}^{2}h_{2}^{3}+h_{2}^{5})(h_{1}k_{1,t}-h_{2}k_{2,t}-k_{1}h_{1,t}+k_{2}h_{2,t})}{(h_{1}^{5}-10h_{1}^{3}h_{2}^{2}+5h_{1}h_{2}^{4})^{2}+(5h_{1}^{4}h_{2}-10h_{1}^{2}h_{2}^{3}+h_{2}^{5})^{2}}(h_{1}k_{1,x}-h_{2}k_{2,x}-k_{1}h_{1,x}+k_{2}h_{2,x})\nonumber\\
+\frac{(h_{1}^{5}-10h_{1}^{3}h_{2}^{2}+5h_{1}h_{2}^{4})(h_{2}k_{1,t}+h_{1}k_{2,t}-k_{2}h_{1,t}-k_{1}h_{2,t})}{(h_{1}^{5}-10h_{1}^{3}h_{2}^{2}+5h_{1}h_{2}^{4})^{2}+(5h_{1}^{4}h_{2}-10h_{1}^{2}h_{2}^{3}+h_{2}^{5})^{2}}(h_{1}k_{1,x}-h_{2}k_{2,x}-k_{1}h_{1,x}+k_{2}h_{2,x})\nonumber\\
+\frac{(h_{1}^{5}-10h_{1}^{3}h_{2}^{2}+5h_{1}h_{2}^{4})(h_{1}k_{1,t}-h_{2}k_{2,t}-k_{1}h_{1,t}+k_{2}h_{2,t})}{(h_{1}^{5}-10h_{1}^{3}h_{2}^{2}+5h_{1}h_{2}^{4})^{2}+(5h_{1}^{4}h_{2}-10h_{1}^{2}h_{2}^{3}+h_{2}^{5})^{2}}(h_{2}k_{1,x}+h_{1}k_{2,x}-k_{2}h_{1,x}-k_{1}h_{2,x})\nonumber\\
+\frac{(5h_{1}^{4}h_{2}-10h_{1}^{2}h_{2}^{3}+h_{2}^{5})(h_{2}k_{1,t}+h_{1}k_{2,t}-k_{1}h_{1,t}-k_{1}h_{2,t})}{(h_{1}^{5}-10h_{1}^{3}h_{2}^{2}+5h_{1}h_{2}^{4})^{2}+(5h_{1}^{4}h_{2}-10h_{1}^{2}h_{2}^{3}+h_{2}^{5})^{2}}(h_{2}k_{1,x}+h_{1}k_{2,x}-k_{2}h_{1,x}-k_{1}h_{2,x}),
\end{eqnarray*}
\begin{eqnarray*}
J_{15}=\frac{(h^{3}_{1}-3h_{1}h^{2}_{2})(k_{1}h_{1,tx}-k_{2}h_{2,tx}+h_{1}k_{1,tx}-h_{2}k_{2,tx}-h_{1,t}k_{1,x}+h_{2,t}k_{2,x}-h_{1,x}k_{1,t}+h_{2,x}k_{2,t})}{(h_{1}^{3}-3h_{1}h_{2}^{2})^2+(3h^{2}_{1}h_{2}-h^{3}_{2})^2}\nonumber\\
+\frac{(3h^{2}_{1}h_{2}-h^{3}_{2})(k_{2}h_{1,tx}+k_{1}h_{2,tx}+h_{2}k_{1,tx}+h_{1}k_{2,tx}-h_{2,t}k_{1,x}-h_{1,t}k_{2,x}-h_{2,x}k_{1,t}-h_{1,x}k_{2,t})}{(h_{1}^{3}-3h_{1}h_{2}^{2})^2+(3h^{2}_{1}h_{2}-h^{3}_{2})^2},
\end{eqnarray*}
\begin{eqnarray*}
J_{16}=\frac{(3h^{2}_{1}h_{2}-h^{3}_{2})(k_{1}h_{1,tx}-k_{2}h_{2,tx}+h_{1}k_{1,tx}-h_{2}k_{2,tx}-h_{1,t}k_{1,x}+h_{2,t}k_{2,x}-h_{1,x}k_{1,t}+h_{2,x}k_{2,t})}{(h_{1}^{3}-3h_{1}h_{2}^{2})^2+(3h^{2}_{1}h_{2}-h^{3}_{2})^2}\nonumber\\
+\frac{(h^{3}_{1}-3h_{1}h^{2}_{2})(k_{2}h_{1,tx}+k_{1}h_{2,tx}+h_{2}k_{1,tx}+h_{1}k_{2,tx}-h_{2,t}k_{1,x}-h_{1,t}k_{2,x}-h_{2,x}k_{1,t}-h_{1,x}k_{2,t})}{(h_{1}^{3}-3h_{1}h_{2}^{2})^2+(3h^{2}_{1}h_{2}-h^{3}_{2})^2},
\end{eqnarray*}
\begin{eqnarray*}
J_{17}=\frac{k_{1}(-6h^{2}_{1}h_{2}^{2}+h^{4}_{1}+h^{4}_{2})+k_{2}(4h^{3}_{1}h_{2}-4h_{1}h^{3}_{2})}{(-6h^{2}_{1}h_{2}^{2}+h^{4}_{1}+h^{4}_{2})^2+(4h^{3}_{1}h_{2}-4h_{1}h^{3}_{2})^2}(h_{1}h_{1,tx}-h_{2}h_{2,tx}-h_{1,t}h_{1,x}+h_{2,t}h_{2,x})\nonumber\\
-\frac{k_{2}(-6h^{2}_{1}h_{2}^{2}+h^{4}_{1}+h^{4}_{2})-k_{1}(4h^{3}_{1}h_{2}-4h_{1}h^{3}_{2})}{(-6h^{2}_{1}h_{2}^{2}+h^{4}_{1}+h^{4}_{2})^2+(4h^{3}_{1}h_{2}-4h_{1}h^{3}_{2})^2}(h_{2}h_{1,tx}+h_{1}h_{2,tx}-h_{2,t}h_{1,x}-h_{1,t}h_{2,x}),
\end{eqnarray*}
\begin{eqnarray*}
J_{18}=\frac{k_{2}(-6h^{2}_{1}h_{2}^{2}+h^{4}_{1}+h^{4}_{2})-k_{1}(4h^{3}_{1}h_{2}-4h_{1}h^{3}_{2})}{(-6h^{2}_{1}h_{2}^{2}+h^{4}_{1}+h^{4}_{2})^2+(4h^{3}_{1}h_{2}-4h_{1}h^{3}_{2})^2}(h_{1}h_{1,tx}-h_{2}h_{2,tx}-h_{1,t}h_{1,x}+h_{2,t}h_{2,x})\nonumber\\
+\frac{k_{1}(-6h^{2}_{1}h_{2}^{2}+h^{4}_{1}+h^{4}_{2})+k_{2}(4h^{3}_{1}h_{2}-4h_{1}h^{3}_{2})}{(-6h^{2}_{1}h_{2}^{2}+h^{4}_{1}+h^{4}_{2})^2+(4h^{3}_{1}h_{2}-4h_{1}h^{3}_{2})^2}(h_{2}h_{1,tx}+h_{1}h_{2,tx}-h_{2,t}h_{1,x}-h_{1,t}h_{2,x}),
\end{eqnarray*}
\begin{eqnarray*}
J_{19}=\frac{\mu_{1}\nu_{1}+\mu_{2}\nu_{2}}{\mu_{1}^{2}+\mu_{2}^{2}}\omega_{1}+\frac{\mu_{2}\nu_{1}-\mu_{1}\nu_{2}}{\mu_{1}^{2}+\mu_{2}^{2}}\omega_{2},\nonumber\\
J_{20}=\frac{\mu_{1}\nu_{2}-\mu_{2}\nu_{1}}{\mu_{1}^{2}+\mu_{2}^{2}}\omega_{1}+\frac{\mu_{1}\nu_{1}+\mu_{2}\nu_{2}}{\mu_{1}^{2}+\mu_{2}^{2}}\omega_{2},\nonumber\\
J_{21}=\frac{\mu_{1}\nu_{3}+\mu_{2}\nu_{4}}{\mu_{1}^{2}+\mu_{2}^{2}}\omega_{3}+\frac{\mu_{2}\nu_{3}-\mu_{1}\nu_{4}}{\mu_{1}^{2}+\mu_{2}^{2}}\omega_{4},\nonumber\\
J_{22}=\frac{\mu_{1}\nu_{4}-\mu_{2}\nu_{3}}{\mu_{1}^{2}+\mu_{2}^{2}}\omega_{3}+\frac{\mu_{1}\nu_{3}+\mu_{2}\nu_{4}}{\mu_{1}^{2}+\mu_{2}^{2}}\omega_{4},
\end{eqnarray*}
where
\begin{eqnarray*}
\mu_{1}=h_{1}^{9}-36h_{1}^{7}h_{2}^{2}+126h_{1}^{5}h_{2}^{4}-84h_{1}^{3}h_{2}^{6}+9h_{1}h_{2}^{8},\nonumber\\
\mu_{2}=9h_{1}^{8}h_{2}-84h_{1}^{6}h_{2}^{3}+126h_{1}^{4}h_{2}^{5}-36h_{1}^{2}h_{2}^{7}+h_{2}^{9},
\end{eqnarray*}
\begin{eqnarray*}
\nu_{1}=k_{2}^{2}h_{2,x}^{2}+2h_{1}k_{2,x}k_{1}h_{2,x}-2h_{2}k_{2,x}k_{2}h_{2,x}+2h_{2}k_{1,x}k_{1}h_{2,x}-4k_{1}h_{1,x}k_{2}h_{2,x}-k_{1}^{2}h_{2,x}^{2}+h_{2}^{2}k_{2,x}^{2}\nonumber\\
+2h_{2}k_{1,x}k_{2}h_{1,x}-2h_{1}k_{1,x}k_{1}h_{1,x}+2h_{2}k_{2,x}k_{1}h_{1,x}+k_{1}^{2}h_{1,x}^{2}+h_{1}^{2}k_{1,x}^{2}+2h_{1}k_{2,x}k_{2}h_{1,x}-h_{2}^{2}k_{1,x}^{2}\nonumber\\
-4h_{1}k_{1,x}h_{2}k_{2,x}-h_{1}^{2}k_{2,x}^{2}-k_{2}^{2}h_{1,x}^{2}+2h_{1}k_{1,x}k_{2}h_{2,x},\nonumber\\
\nu_{2}=-2k_{2}h_{2,x}^{2}k_{1}-2h_{1}k_{1,x}k_{1}h_{2,x}-2k_{1}h_{1,x}h_{1}k_{2,x}+2h_{2}k_{2,x}k_{2}h_{1,x}-2k_{2}^{2}h_{2,x}h_{1,x}+2h_{2}k_{2,x}k_{1}h_{2,x}\nonumber\\
+2k_{2}h_{2,x}h_{2}k_{1,x}-2h_{2}^{2}k_{2,x}k_{1,x}+2k_{2}h_{2,x}h_{1}k_{2,x}+2h_{1}k_{1,x}^{2}h_{2}-2h_{1}k_{1,x}k_{2}h_{1,x}-2h_{2}k_{2,x}^{2}h_{1}\nonumber\\
+2h_{1}^{2}k_{1,x}k_{2,x}+2k_{1}h_{1,x}^{2}k_{2}-2k_{1}h_{1,x}h_{2}k_{1,x}+2k_{1}^{2}h_{1,x}h_{2,x},\nonumber\\
\nu_{3}=-2h_{2}k_{2,t}k_{2}h_{2,t}-k_{1}^{2}h_{2,t}^{2}+2h_{2}k_{1,t}k_{2}h_{1,t}-h_{2}^{2}k_{1,t}^{2}+2h_{2}k_{2,t}k_{1}h_{1,t}+h_{1}^{2}k_{1,t}^{2}-2h_{1}k_{1,t}k_{1}h_{1,t}\nonumber\\
+k_{1}^{2}h_{1,t}^{2}+k_{2}^{2}h_{2,t}^{2}-4h_{1}k_{2,t}h_{2}k_{1,t}-k_{2}^{2}h_{1,t}^{2}+2h_{1}k_{1,t}k_{2}h_{2,t}+2h_{1}k_{2,t}k_{1}h_{2,t}-h_{1}^{2}k_{2,t}^{2}\nonumber\\
-4k_{1}h_{1,t}k_{2}h_{2,t}+2h_{2}k_{1,t}k_{1}h_{2,t}+2h_{1}k_{2,t}k_{2}h_{1,t}+h_{2}^{2}k_{2,t}^{2},\nonumber\\
\nu_{4}=2h_{2}k_{2,t}k_{2}h_{1,t}+2h_{1}k_{1,t}^{2}h_{2}+2h_{1}k_{2,t}k_{2}h_{2,t}+2k_{1}h_{1,t}^{2}k_{2}-2h_{1}k_{2,t}k_{1}h_{1,t}+2h_{1}^{2}k_{1,t}k_{2,t}\nonumber\\
-2k_{1}h_{2,t}^{2}k_{2}-2h_{2}^{2}k_{1,t}k_{2,t}+2h_{2}k_{2,t}k_{1}h_{2,t}+2h_{2}k_{1,t}k_{2}h_{2,t}+2k_{1}^{2}h_{1,t}h_{2,t}-2k_{2}^{2}h_{1,t}h_{2,t}\nonumber\\
-2h_{2}k_{1,t}k_{1}h_{1,t}-2h_{1}k_{1,t}k_{1}h_{2,t}-2h_{1}k_{1,t}k_{2}h_{1,t}-2h_{1}k_{2,t}^{2}h_{2},
\end{eqnarray*}
and
\begin{eqnarray*}
\omega_{1}=(h_{1}k_{1}-h_{2}k_{2})h_{1,tt}-(h_{2}k_{1}+h_{1}k_{2})h_{2,tt}+(-h_{1}^{2}+h_{2}^{2})k_{1,tt}+2h_{1}h_{2}k_{2,tt}-3k_{1}(h_{1,t}^{2}-h_{2,t}^{2})\nonumber\\
+6k_{2}h_{1,t}h_{2,t}+(3h_{1}h_{1,t}-3h_{2}h_{2,t})k_{1,t}-(3h_{2}h_{1,t}+3h_{1}h_{2,t})k_{2,t},\nonumber\\
\omega_{2}=(h_{2}k_{1}+h_{1}k_{2})h_{1,tt}+(h_{1}k_{1}-h_{2}k_{2})h_{2,tt}+(-h_{1}^{2}+h_{2}^{2})k_{2,tt}-2h_{1}h_{2}k_{1,tt}-3k_{2}(h_{1,t}^{2}-h_{2,t}^{2})\nonumber\\
-6k_{1}h_{1,t}h_{2,t}+(3h_{2}h_{1,t}+3h_{1}h_{2,t})k_{1,t}+(3h_{1}h_{1,t}-3h_{2}h_{2,t})k_{2,t},\nonumber\\
\omega_{3}=(h_{1}k_{1}-h_{2}k_{2})h_{1,xx}-(h_{2}k_{1}+h_{1}k_{2})h_{2,xx}+(-h_{1}^{2}+h_{2}^{2})k_{1,xx}+2h_{1}h_{2}k_{2,xx}-3k_{1}(h_{1,x}^{2}-h_{2,x}^{2})\nonumber\\
+6k_{2}h_{1,x}h_{2,x}+(3h_{1}h_{1,x}-3h_{2}h_{2,x})k_{1,x}-(3h_{2}h_{1,x}+3h_{1}h_{2,x})k_{2,x},\nonumber\\
\omega_{4}=(h_{2}k_{1}+h_{1}k_{2})h_{1,xx}+(h_{1}k_{1}-h_{2}k_{2})h_{2,xx}+(-h_{1}^{2}+h_{2}^{2})k_{2,xx}-2h_{1}h_{2}k_{1,xx}-3k_{2}(h_{1,x}^{2}-h_{2,x}^{2})\nonumber\\
-6k_{1}h_{1,x}h_{2,x}+(3h_{2}h_{1,x}+3h_{1}h_{2,x})k_{1,x}+(3h_{1}h_{1,x}-3h_{2}h_{2,x})k_{2,x},
\end{eqnarray*}
which are found to be associated with the system of two linear
hyperbolic-type PDEs (\ref{chypsysgen}). These also can be observed to be the complex split of the joint invariants (16).

\section{Applications}
In this section a few examples of systems of hyperbolic-type
equations are provided to illustrate the invariance criteria
developed.

\textbf{1.} A system of two hyperbolic-type PDEs
\begin{eqnarray}
u_{tx}+\left(a_{1}-\frac{1}{x}\right)u_{t}-a_{2}v_{t}+\left(b_{1}+\frac{2}{t}\right)u_{x}-b_{2}v_{x}+\left(c_{1}-\frac{b_{1}}{x}+2\frac{a_{1}}{t}-\frac{2}{tx}\right)u\nonumber\\
-\left(c_{2}-\frac{b_{2}}{x}+2\frac{a_{2}}{t}\right)v=0,~~~~~~~~~~~~~~~~~~~~~~~~~~~~~~~~~~~~~~~~~~~~~~~~~~~~~~~~\nonumber\\
v_{tx}+a_{2}u_{t}+\left(a_{1}-\frac{1}{x}\right)v_{t}+b_{2}u_{x}+\left(b_{1}+\frac{2}{t}\right)v_{x}+\left(c_{2}-\frac{b_{2}}{x}+2\frac{a_{2}}{t}\right)u~~~~~~~~\nonumber\\
+\left(c_{1}-\frac{b_{1}}{x}+2\frac{a_{1}}{t}-\frac{2}{tx}\right)v=0,~~~~~~~~~~~~~~~~~~~~~~~~~~~~~~~~~~~~~~~~~~~~~~~~~\label{ex1hsys}
\end{eqnarray}
corresponds to a complex hyperbolic equation in two independent
variables
\begin{eqnarray}
w_{tx}+\left(a-\frac{1}{x}\right)w_{t}+\left(b+\frac{2}{t}\right)w_{x}+\left(c-\frac{b}{x}+2\frac{a}{t}-\frac{2}{tx}\right)w=0,\label{ex1hcom}
\end{eqnarray}
where $a$ is a complex constant $a=a_{1}+ia_{2}$. The following
complex transformation of the dependent variable
$w=(x/t^2)\overline{w}$ maps the complex equation (\ref{ex1hcom}) to
\begin{eqnarray}
\overline{w}_{tx}+a\overline{w}_{t}+b\overline{w}_{x}+c\overline{w}=0.\label{ex1trhcom}
\end{eqnarray}
Both the complex hyperbolic equations (\ref{ex1hcom}) and
(\ref{ex1trhcom}) are transformable to each other because they have the
same semi-invariants
\begin{eqnarray}
h=ab-c=k.
\end{eqnarray}
The system of hyperbolic-type equations (\ref{ex1hsys}) is
transformable to
\begin{eqnarray}
\overline{u}_{tx}+a_{1}\overline{u}_{t}-a_{2}\overline{v}_{t}+b_{1}\overline{u}_{x}-b_{2}\overline{v}_{x}+c_{1}\overline{u}-c_{2}\overline{v}=0,\nonumber\\
\overline{v}_{tx}+a_{2}\overline{u}_{t}+a_{1}\overline{v}_{t}+b_{2}\overline{u}_{x}+b_{1}\overline{v}_{x}+c_{2}\overline{u}+c_{1}\overline{v}=0,\label{ex1trhsys}.
\end{eqnarray}
The real transformations of the dependent variables
\begin{eqnarray}
u=(x/t^2)\overline{u},~~v=(x/t^2)\overline{v},
\end{eqnarray}
are obtained by splitting the complex dependent transformation used
to map the complex equations (\ref{ex1hcom}) and (\ref{ex1trhcom})
into each other. Semi-invariants associated with both the systems
(\ref{ex1trhsys}) are
\begin{eqnarray}
h_{1}=a_{1}b_{1}-a_{2}b_{2}-c_{1}=k_{1},\nonumber\\
h_{2}=a_{1}b_{2}+a_{2}b_{1}-c_{2}=k_{2},
\end{eqnarray}
which guarantees that both the systems are mappable into each other.

\textbf{2.} An uncoupled system of PDEs
\begin{eqnarray}
u_{z_{1}z_{2}}+2az_{1}^{2}u_{z_{1}}+2bz_{1}u_{z_{2}}+4cz_{1}u=0,\nonumber\\
v_{z_{1}z_{2}}+2az_{1}^{2}v_{z_{1}}+2bz_{1}v_{z_{2}}+4cz_{1}v=0,\label{ex2unsys}
\end{eqnarray}
is transformable to
\begin{eqnarray}
u_{tx}+atu_{t}+bu_{x}+cu=0,\nonumber\\
v_{tx}+atv_{t}+bv_{x}+cv=0,\label{ex2untrsys}
\end{eqnarray}
via invertible transformations of the independent variables
\begin{eqnarray}
z_{1}=\sqrt{t},~~z_{2}=\frac{1}{2}(x-1).\label{ex2trans}
\end{eqnarray}
These are the invertible maps that also reduce the base complex
hyperbolic equation of the form
\begin{eqnarray}
w_{z_{1}z_{2}}+2az_{1}^{2}w_{z_{1}}+2bz_{1}w_{z_{2}}+4cz_{1}w=0,\label{ex2hcom}
\end{eqnarray}
with the semi-invariants
\begin{eqnarray}
I_{1}=\frac{c}{abz_{1}^{2}},~~I_{2}=bz_{1}^{2},~~I_{3}=0,~~I_{4}=\frac{c}{a},~~I_{5}=0=I_{6},\label{ex2sicom}
\end{eqnarray}
to a simple linear form
\begin{eqnarray}
w_{tx}+atw_{t}+bw_{x}+cw=0,
\end{eqnarray}
with the following semi-invariants
\begin{eqnarray}
I_{1}=\frac{c}{abt},~~I_{2}=bt,~~I_{3}=0,~~I_{4}=\frac{c}{a},~~I_{5}=0=I_{6}.\label{ex2sitrcom}
\end{eqnarray}
Notice that the semi-invariants (\ref{ex2sicom}) and
(\ref{ex2sitrcom}) are the same by means of the transformations of
the independent variables (\ref{ex2trans}). The complex hyperbolic
equation (\ref{ex2hcom}) does not only yield an uncoupled system of
the hyperbolic-type equations (\ref{ex2unsys}). In fact it gives a
coupled system
\begin{eqnarray}
u_{z_{1}z_{2}}+2a_{1}z_{1}^{2}u_{z_{1}}-2a_{2}z_{1}^{2}v_{z_{1}}+2b_{1}z_{1}u_{z_{2}}-2b_{2}z_{1}v_{z_{2}}+4c_{1}z_{1}u-4c_{2}z_{1}v=0,\nonumber\\
v_{z_{1}z_{2}}+2a_{2}z_{1}^{2}u_{z_{1}}+2a_{1}z_{1}^{2}v_{z_{1}}+2b_{2}z_{1}u_{z_{2}}+2b_{1}z_{1}v_{z_{2}}+4c_{2}z_{1}u+4c_{1}z_{1}v=0.\label{ex2cupsys}
\end{eqnarray}
This system of two hyperbolic-type equations can be mapped to
\begin{eqnarray}
u_{tx}+a_{1}tu_{t}-a_{2}tv_{t}+b_{1}u_{x}-b_{2}v_{x}+c_{1}u-c_{2}v=0,\nonumber\\
v_{tx}+a_{2}tu_{t}+a_{1}tv_{t}+b_{2}u_{x}+b_{1}v_{x}+c_{2}u+c_{1}v=0,\label{ex2cuptrsys}
\end{eqnarray}
under the transformations (\ref{ex2trans}) that are already used to
map the base complex equation to its canonical form.

\textbf{3.} Invoking the following change of the independent
variables
\begin{eqnarray}
z_{1}=e^{t},~~z_{2}=\sqrt{x},\label{ex3trans}
\end{eqnarray}
in a coupled system of two hyperbolic-type equations of the form
\begin{eqnarray}
u_{z_{1},z_{2}}+2a_{1}z_{2}\ln{z_{1}}u_{z_{1}}-2a_{2}z_{2}\ln{z_{1}}v_{z_{1}}+\frac{b_{1}}{z_{1}}u_{z_{2}}-\frac{b_{2}}{z_{1}}v_{z_{2}}+\frac{2c_{1}z_{2}}{z_{1}}u-\frac{2c_{2}z_{2}}{z_{1}}v=0,\nonumber\\
v_{z_{1},z_{2}}+2a_{2}z_{2}\ln{z_{1}}u_{z_{1}}+2a_{1}z_{2}\ln{z_{1}}v_{z_{1}}+\frac{b_{2}}{z_{1}}u_{z_{2}}+\frac{b_{1}}{z_{1}}v_{z_{2}}+\frac{2c_{2}z_{2}}{z_{1}}u+\frac{2c_{1}z_{2}}{z_{1}}v=0,
\end{eqnarray}
transforms it to
\begin{eqnarray}
u_{tx}+a_{1}tu_{t}-a_{2}tv_{t}+b_{1}u_{x}-b_{2}v_{x}+c_{1}u-c_{2}v=0,\nonumber\\
v_{tx}+a_{2}tu_{t}+a_{1}tv_{t}+b_{2}u_{x}+b_{1}v_{x}+c_{2}u+c_{1}v=0.
\end{eqnarray}
The transformation of these systems under the invertible change of
the independent variables follows from the base complex hyperbolic
equation
\begin{eqnarray}
w_{z_{1}z_{2}}+2a\ln{z_{1}}w_{z_{1}}+\frac{b}{z_{1}}w_{z_{2}}+\frac{2cz_{2}}{z_{1}}w=0.
\end{eqnarray}
It can be transformed to another linear form
\begin{eqnarray}
w_{tx}+atw_{t}+bw_{x}+cw=0,
\end{eqnarray}
under the invertible transformations (\ref{ex3trans}). Similarly,
the invertible transformations of the independent variables
(\ref{ex3trans}) map the following system of PDEs
\begin{eqnarray}
u_{z_{1},z_{2}}+2a_{1}z_{2}u_{z_{1}}-2a_{2}z_{2}v_{z_{1}}+\frac{b_{1}}{z_{1}}u_{z_{2}}-\frac{b_{2}}{z_{1}}v_{z_{2}}+\frac{2c_{1}z_{2}}{z_{1}}u-\frac{2c_{2}z_{2}}{z_{1}}v=0,\nonumber\\
v_{z_{1},z_{2}}+2a_{2}z_{2}u_{z_{1}}+2a_{1}z_{2}v_{z_{1}}+\frac{b_{2}}{z_{1}}u_{z_{2}}+\frac{b_{1}}{z_{1}}v_{z_{2}}+\frac{2c_{2}z_{2}}{z_{1}}u+\frac{2c_{1}z_{2}}{z_{1}}v=0,
\end{eqnarray}
to
\begin{eqnarray}
u_{tx}+a_{1}u_{t}-a_{2}v_{t}+b_{1}u_{x}-b_{2}v_{x}+c_{1}u-c_{2}v=0,\nonumber\\
v_{tx}+a_{2}u_{t}+a_{1}v_{t}+b_{2}u_{x}+b_{1}v_{x}+c_{2}u+c_{1}v=0.
\end{eqnarray}

\textbf{4.} Consider an uncoupled system of two hyperbolic type PDEs
\begin{eqnarray}
g_{1,tx}+\frac{\lambda}{2}(g_{1,t}+g_{1,x})=0,\nonumber\\
g_{2,tx}+\frac{\lambda}{2}(g_{2,t}+g_{2,x})=0,\label{ex1sysun}
\end{eqnarray}
for which $h_{1}=k_{1}=\frac{\lambda^2}{4}$, and $h_{2}=k_{2}=0$. This
implies that
\begin{eqnarray}
J_{1}=1,~~J_{2}=\cdots=J_{12}=0.\label{ex1inv}
\end{eqnarray}
The system (\ref{ex1sysun}) is transformable to another system with
the same invariants as given in (\ref{ex1inv}) where
$h_{1}=k_{1}=-1$, $h_{2}=k_{2}=0$. The transformed system reads as
\begin{eqnarray}
f_{1,z_{1}z_{2}}+f_{1}=0,\nonumber\\
f_{2,z_{1}z_{2}}+f_{2}=0.\label{ex1sysuntr}
\end{eqnarray}
The correspondence between the systems (\ref{ex1sysun}) and
(\ref{ex1sysuntr}) is established by
\begin{eqnarray}
z_{1}=\frac{\lambda}{2}t,~~z_{2}=-\frac{\lambda}{2}x,~~f_{1}=g_{1}\exp(\frac{\lambda
t+\lambda x}{2}),~~f_{2}=g_{2}\exp(\frac{\lambda t+\lambda x}{2}).
\end{eqnarray}
These transformations are obtainable from
\begin{eqnarray}
z_{1}=\frac{\lambda}{2}t,~~z_{2}=-\frac{\lambda}{2}x,~~w=u\exp(\frac{\lambda
t+\lambda x}{2}),\label{ex1comtran}
\end{eqnarray}
by  $w=f_{1}+if_{2}$ and $u=g_{1}+ig_{2}$. The complex
transformations map the complex scalar PDE
\begin{eqnarray}
w_{,z_{1}z_{2}}+\frac{\lambda}{2}(w_{,z_{1}}+w_{,z_{2}})=0,\label{ex1sca}
\end{eqnarray}
with $h=k=\frac{\lambda^2}{4}$ and $p=1$, to an equation
\begin{eqnarray}
u_{,tx}+u=0,
\end{eqnarray}
for which $h=k=-1$ and $p=1$. Notice that the substitution
$\lambda=\lambda_{1}+i\lambda_{2}$, in the equation (\ref{ex1sca})
results in a coupled system of two hyperbolic-type PDEs but it can
not be transformed by the complex method. The reason is the complex
transformations (\ref{ex1comtran}) where the two independent
variables split into four add extra dimensions. Therefore, the
complex procedure fails for that case.

\textbf{5.} The complex transformations of the form
\begin{eqnarray}
z_{1}=\frac{1}{t},~~z_{2}=2x,~~w=\frac{u}{x},\label{ex2comtran}
\end{eqnarray}
map the following Lie canonical form
\begin{eqnarray}
w_{,z_{1}z_{2}}+\alpha z^{2}_{2}w_{,z_{2}}+2w=0,\label{ex2trsca}
\end{eqnarray}
to
\begin{eqnarray}
u_{,tx}-\frac{1}{x}u_{,t}-\frac{\alpha
x^2}{t^2}u_{,x}+\frac{1}{t^2}(\alpha x-2)u=0.\label{ex2sca}
\end{eqnarray}
The invariant quantities associated with both the scalar Lie
canonical form and the hyperbolic equations are $h=-1,~k=2\alpha
x-1,~p=2(1-\alpha x)$ and $h=2/t^2,~k=\frac{2(1-\alpha
x)}{t^2},~p=1-\alpha x$, respectively. Inserting $u=g_{1}+ig_{2}$ in
the equation (\ref{ex2sca}) while keeping $\alpha$ \emph{a real
constant} yields an uncoupled system of two PDEs
\begin{eqnarray}
g_{1,tx}-\frac{1}{x}g_{1,t}-\frac{\alpha
x^2}{t^2}g_{1,x}+\frac{\alpha x-2}{t^2}g_{1}=0,\nonumber\\
g_{2,tx}-\frac{1}{x}g_{2,t}-\frac{\alpha
x^2}{t^2}g_{2,x}+\frac{\alpha x-2}{t^2}g_{2}=0.\label{ex2sys}
\end{eqnarray}
The system (\ref{ex2sys}) is transformable to another system of the
form
\begin{eqnarray}
f_{1,z_{1}z_{2}}+\alpha x^2f_{1,z_{2}}+2f_{1}=0,\nonumber\\
f_{2,z_{1}z_{2}}+\alpha x^2f_{2,z_{2}}+2f_{2}=0,\label{ex2trsys}
\end{eqnarray}
under a change of the dependent and independent variables
\begin{eqnarray}
z_{1}=\frac{1}{t},~~z_{2}=2x,~~f_{1}=\frac{g_{1}}{x},~~f_{2}=\frac{g_{2}}{x}.\label{ex2retrans}
\end{eqnarray}
These transformations are the real and imaginary parts of the
complex transformations (\ref{ex2comtran}) and the transformed
system is obtained by splitting the Lie canonical form
(\ref{ex2trsca}) into the real and imaginary parts. The invariance
criteria that ensure such a transformation of the system are
satisfied. These quantities for both the systems (\ref{ex2sys}) and
(\ref{ex2trsys}) are
\begin{eqnarray}
h_{1}=\frac{2}{t^2},~~k_{1}=\frac{2(1-\alpha
x)}{t^2},~~h_{2}=0=k_{2},~~p=\frac{-1}{\alpha x-1},
\end{eqnarray}
and
\begin{eqnarray}
h_{1}=-2,~~k_{1}=2(\alpha x-1),~~h_{2}=0=k_{2},~~p=\frac{1}{1-\alpha
x},
\end{eqnarray}
respectively.

A coupled system
\begin{eqnarray}
g_{1,tx}-\frac{1}{x}g_{1,t}-\frac{\alpha_{1}
x^2}{t^2}g_{1,x}+\frac{\alpha_{2}
x^2}{t^2}g_{2,x}+\frac{\alpha_{1} x-2}{t^2}g_{1}-\frac{\alpha_{2} x}{t^2}g_{2}=0,\nonumber\\
g_{2,tx}-\frac{1}{x}g_{2,t}-\frac{\alpha_{2}
x^2}{t^2}g_{1,x}-\frac{\alpha_{1} x^2}{t^2}g_{2,x}+\frac{\alpha_{2}
x}{t^2}g_{1}+\frac{\alpha_{1} x-2}{t^2}g_{2}=0,\label{ex2cupsys}
\end{eqnarray}
with the invariants
\begin{eqnarray}
h_{1}=\frac{2}{t^2},~~k_{1}=\frac{2(1-\alpha_{1}
x)}{t^2},~~h_{2}=0,~~k_{2}=\frac{-2\alpha_{2}
x}{t^2},\nonumber\\
J_{1}=\frac{1-\alpha_{1}x}{(1-\alpha_{1}x)^{2}+\alpha^{2}_{2}x^2},~~J_{2}=\frac{\alpha_{2}x}{(1-\alpha_{1}x)^{2}+\alpha^{2}_{2}x^2},
\end{eqnarray}
is obtainable from the complex scalar PDE (\ref{ex2sca}) when
$\alpha$ is also complex, i.e., $\alpha=\alpha_{1}+i\alpha_{2}$.
Employing the transformations (\ref{ex2retrans}) on
(\ref{ex2cupsys}) one arrives at a coupled system
\begin{eqnarray}
f_{1,z_{1}z_{2}}+\alpha_{1} x^2f_{1,z_{2}}-\alpha_{2} x^2f_{2,z_{2}}+2f_{1}=0,\nonumber\\
f_{2,z_{1}z_{2}}+\alpha_{2}x^2f_{1,z_{2}}+\alpha_{1}x^2f_{2,z_{2}}+2f_{2}=0,\label{ex2trcupsys}
\end{eqnarray}
which is the real analogue of the complex transformed equation
(\ref{ex2trsca}) and satisfies the invariance criteria, where
\begin{eqnarray}
h_{1}=-2,~~k_{1}=2(\alpha_{1} x-1),~~h_{2}=0,~~k_{2}=2\alpha_{2}
x,\nonumber\\
J_{1}=\frac{1-\alpha_{1}x}{(1-\alpha_{1}x)^{2}+\alpha^{2}_{2}x^2},~~J_{2}=\frac{\alpha_{2}x}{(1-\alpha_{1}x)^{2}+\alpha^{2}_{2}x^2}.
\end{eqnarray}

\section{Conclusion}
Semi-invariants of the hyperbolic and parabolic PDEs in
two independent variables have been obtained by transforming the dependent
or independent variables \cite{book1,book2,ibrlp,lap,fmjhs}. Further, the infinitesimal approach has
been utilized to derive the joint invariants for the hyperbolic
and parabolic equations \cite{fmjhs,ibrc,mah,jf,ibrb}. The semi-invariants of the hyperbolic and
parabolic PDEs have been extended to systems of such equations by
complex symmetry analysis \cite{comh,psysp}. The real and complex approaches were
investigated in this work for the invariants of a system of two linear hyperbolic
equations.

Semi-invariants of a special class of systems of two hyperbolic-type
PDEs were derived here using real and complex methods developed for
such systems of equations. Both the procedures are adopted to find
the semi-invariants of the system of two hyperbolic-type equations
that is obtainable from a complex hyperbolic PDE. Semi-invariants
associated with the invertible change of the dependent as well as
independent variables are deduced by both the real and complex
methods. It is shown that same invariant quantities for the system
of hyperbolic-type PDEs appear due to complex and real procedures,
in the case of transformations of only the dependent variables.
However, the semi-invariants of this system obtained by real
symmetry analysis are different from those provided by the complex procedure.
Furthermore, the joint invariants of this system of hyperbolic-type
equations obtained by both the methods are also found to be different.

{\bf Acknowledgments}. FM is thankful to the NRF of South Africa for  an
enabling research grant. AA thanks DECMA of Wits and NUST for support during the time
this work was completed.

\end{document}